\documentclass[12pt]{article}        
\usepackage[french]{babel}     
\newtheorem{lem}{Lemme}[section]
\newtheorem{prop}{Proposition}[section]
\newtheorem{thm}{Th\'eor\`eme}[section]
\newtheorem{cor}{Corollaire}[section]
\newtheorem{conj}{Conjecture}[section]
{\itshape}{\rm}
\newtheorem{rem}{Remarque}[section]
\newtheorem{defi}{D\'efinition}[section]

\let\fd=\rightarrow
\let\lfd=\longrightarrow
\let\lfc=\longmapsto

\let\Om=\Omega

\let\eps=\epsilon

\let\scr=\scriptstyle
\let\dps=\displaystyle

\let\bs=\bigskip
\let\bb=\bigbreak

\def\N{{\msb N}}
\def\Z{{\msb Z}}

\def\R{{\msb R}}
\def\C{{\msb C}}

\def\F*2g{{\msb F}^*_{2^g}}
\def\f#1{{\msb F}_{#1}}

\def\F#1{{\msb F}_{2^{#1}}}

\def\pt#1{\left(#1\right)}
\def\ssi{si et seulement si }

\def\sg#1{\hbox{$\underline{\hbox{#1}}$}}

\def\biq#1{\quad \hbox{ #1 }\quad }

\def\ie{c'est-\`a-dire }

\newcount\annee\annee=\year\advance\annee by -1900

\font\tenmsb=msbm10 
\font\sevenmsb=msbm7
\font\fivemsb=msbm5
\newfam\msbfam
\def\msb{\fam\msbfam\tenmsb}%
\textfont\msbfam=\tenmsb \scriptfont\msbfam=\sevenmsb%
\scriptscriptfont\msbfam=\fivemsb%

\font\tenmsb=msbm10 at 12pt
\font\sevenmsb=msbm9
\font\fivemsb=msbm6
\newfam\msbfam%
\def\msb{\fam\msbfam\tenmsb}%
\textfont\msbfam=\tenmsb \scriptfont\msbfam=\sevenmsb%
\scriptscriptfont\msbfam=\fivemsb%

\def\ch{\widehat}
\def\E{{\cal E}}

\begin{document}

\title{Sur la non-lin\'earit\'e des fonctions bool\'eennes}
\author{Fran\c cois Rodier\\Institut de Math\'ematiques de Luminy --
C.N.R.S.\\Marseille -- France}
\date{}
\maketitle

\section{Introduction}

Les fonctions bool\'eennes sur l'espace $\f2^m$ interviennent aussi bien
dans la th\'eorie des codes correcteurs d'erreurs (par exemple dans les
codes de Reed-Muller) qu'en cryptographie pour r\'ealiser des syst\`emes de
chiffrement \`a clef secr\`ete.

Dans ces deux cas, les propri\'et\'es des syst\`emes ainsi construits d\'ependent
en particulier de la non-lin\'earit\'e d'une fonction bool\'eenne, concept que je
d\'efinirai  pr\'ecis\'ement plus loin (\S\
\ref{defnl}).
La non-lin\'earit\'e est  li\'ee au rayon de recouvrement des codes de
Reed-Muller.
 C'est aussi un  param\`etre cryptographique important: dans leur article
\cite{cv}, F.~Chabaud et
S.~Vaudenay montrent  que la non-lin\'earit\'e est
un crit\`ere important de r\'esistance aux attaques diff\'erentielles et
lin\'eaires; dans sa th\`ese,  C.~Fontaine \cite{cf}   met en valeur 
l'importance de la non-lin\'earit\'e en cryptographie pour plusieurs syst\`emes
de chiffrement.

Il est utile de pouvoir disposer de fonctions
bool\'eennes ayant la plus grande non-lin\'earit\'e possible,
comme l'ont montr\'e W. Meier et O. Staffelbach dans
\cite{ms}, ainsi que K. Nyberg dans \cite{ny}. Ces fonctions
ont \'et\'e \'etudi\'ees dans le cas o\`u $m$ est pair, sous le nom de fonctions
``courbes'' (cf. J.~Dillon \cite{di}). Leur degr\'e de non-lin\'earit\'e est alors
bien connu, on sait  construire plusieurs s\'eries de fonctions
courbes, mais on ne conna\^it pas encore ni leur
nombre, ni leur classification (cf. les travaux de C.~Carlet, en particulier
l'article de  C.~Carlet et P. Guillot \cite{cg}, ou de  C. Carlet et
A. Klapper \cite{ck}).
Dans le cas o\`u $m$ est impair, la situation est bien diff\'erente: on ne
conna\^it alors la valeur de la non-lin\'earit\'e maximale que pour quelques
valeurs de $m$, et on n'a qu'une conjecture pour les autres valeurs
 (voir le m\'emoire d'habilitation de P.~Langevin \cite{la}).

Dans cet article, je veux montrer que pour traiter les fonctions
bool\'een\-nes, on peut s'inspirer de la th\'eorie des polyn\^omes al\'eatoires qui
a \'et\'e un sujet d'\'etude depuis les travaux de Paley et Zygmund. En effet, le
probl\`eme de la recherche du maximum du degr\'e de non-lin\'earit\'e revient \`a
minimiser la transform\'ee de Fourier de fonctions bool\'eennes. C'est un
probl\`eme analogue aux s\'eries de Fourier sur un tore, o\`u l'on cherche \`a
minimiser la transform\'ee de Fourier des fonctions sur $\Z$ prenant les
valeurs
$\pm1$ pour un ensemble fini (et 0 ailleurs) , ce qui revient \`a chercher
 \`a mini\-miser les valeurs  des polyn\^omes \`a coefficient $\pm1$
(polyn\^omes {\sl al\'eatoires}) sur l'ensemble des nombres complexes de module
1.

Dans cet article, on s'inspire des travaux de R.~Salem et A.~Zygmund
\cite{sz} et de J-P.~Kahane \cite{ka} sur les polyn\^omes al\'eatoires, en les
transposant sur les fonctions bool\'eennes. On trouve ainsi une \'evaluation de
la moyenne des normes dans
$L_\infty$ des transform\'ees de Fourier des fonctions bool\'eennes,
qui n'est pas trop \'eloign\'ee de sa valeur minimale th\'eorique, $2^{m/2}$.
Cela donne une
\'evaluation  de la moyenne des degr\'es de non-lin\'earit\'e de ces fonctions.
On retrouve en particulier  le fait  que la plupart des
fonctions bool\'eennes ont une grande non-lin\'earit\'e, un
r\'esultat mis en \'evidence r\'ecemment par D.~Olej\'ar et M.~Stanek \cite{os} et
C. Carlet
\cite{ca2, ca3} (cf. th\'eor\`eme \ref{th1}). Le r\'esultat que 
j'ai d\'emontr\'e implique en outre que presque toutes les fonctions bool\'eennes
ont un non-lin\'earit\'e voisine d'une m\^eme valeur. Cette propri\'et\'e est
illustr\'ee par exemple par les diagrammes de
\cite{am}, qui exhibent la non-lin\'earit\'e de fonctions bool\'eennes en vue de
la construction de bo\^ites de substitutions (s-boxes) utilis\'ees dans les
chiffrements par blocs, ou de l'\'etude statistique de \cite{cf} chapitre 6.

De plus, en transposant une \'etude de D. Newman et J. Byrnes \cite{nb} sur
les normes dans
$L_4$ des polyn\^omes,  nous avons \'et\'e amen\'es \`a \'etudier une conjecture sur
la norme dans $L_4$ des transform\'ee de Fourier de fonctions
bool\'eennes.
On retrouve ainsi le crit\`ere de la ``somme des carr\'es'', reli\'e au
crit\`ere de propagation, pour les fonctions bool\'eennes.  
Ce crit\`ere a \'et\'e \'etudi\'e par Xian-Mo Zhang et Yuliang Zheng \cite{zz}, ou par
P. St\u anic\u a \cite{st}. Son rapport avec la non-lin\'earit\'e a \'et\'e \'etudi\'e
par A.~Canteaut et al. \cite{cccf}.

\section{Pr\'eliminaires}

\subsection{Fonctions bool\'eennes}

Soit $m$ un entier positif et $q=2^m$.

\begin{defi}
Une fonction bool\'eenne \`a $m$ variables est une application de l'espace
$V_m=(\f2)^m$ dans $\f2$.
\end{defi}

Une fonction bool\'eenne est \emph{lin\'eaire} si c'est une forme lin\'eaire sur
l'espace vectoriel $(\f2)^m$. Elle est dite \emph{affine} si elle est
\'egale \`a une fonction lin\'eaire \`a une constante pr\`es. 

\subsection{Rayon
de recouvrement du code du Reed-Muller du premier ordre et amplitude
spectrale}

\begin{defi}
L'amplitude spectrale de la fonction bool\'eenne $g$ est \'egale \`a
$$  S(g) = 
\sup_{v\in V_m} \Bigl| \sum_{x\in V_m}(-1)^{(g(x)+v\cdot x)}\Bigr|$$
o\`u $v\cdot x$ note le produit scalaire usuel dans $V_m$.
C'est le maximum de la transform\'ee de Fourier de  $(-1)^g$.
\end{defi}

Cette amplitude spectrale est reli\'ee au rayon
de recouvrement du code du Reed-Muller.

En effet, un code de  Reed et Muller ${\cal R}_m $ d'ordre 1 sur $V_m$ est
l'espace vectoriel  des fonctions bool\'eennes affines sur $V_m$. 
  Le   rayon
de recouvrement $r_m$ du code est le plus petit entier tel que chaque
vecteur de longueur $2^m$ (c'est \`a dire chaque fonction $V_m\fd\f2$) est \`a une 
distance (de Hamming) d'un mot de code de ${\cal R}_ m $ au plus \'egale \`a 
$r_m$.
  On v\'erifie que
  $$r_m = 2^{m-1} - {1\over 2}\mu _m 
\biq{o\`u}
\mu _m = \inf_g  S(g) 
$$
  o\`u $g$ est une fonction $V_m \fd \f 2$
et o\`u 
$ S(g)$ est l'amplitude spectrale de la fonction $g$.

\subsection{Non-lin\'earit\'e}
\label{defnl}

\begin{defi}
On appelle degr\'e de non-lin\'earit\'e d'une fonction bool\'eenne $g$ \`a
$m$ variables et on le note $nl(g)$ la distance qui la s\'epare de
l'ensemble des fonctions affines \`a $m$ variables :
$$nl(g) = \min_{h \hbox{\,\scriptsize affine }} d(g, h)$$
o\`u $d$ est la distance de Hamming.
\end{defi}

  \begin{prop}
Soit $g$ une fonction bool\'eenne \`a m variables. Son degr\'e de
non-lin\'earit\'e est \'egal \`a
$$nl(g)  = 2^{m-1} - {1\over 2}S(g).
$$
\end{prop}

{\sl D\'emonstration} --

C'est la m\^eme d\'emonstration que pour le rayon de recouvrement d'un code
de Reed et Muller.

\subsection{R\'esultats connus, conjecture}

Le rayon
de recouvrement du code du Reed-Muller du premier ordre est bien connu
pour une dimension $m$ paire: 
$ \mu _m $ vaut $2^{m/2}$.
Pour $m$ impair, on n'a connu longtemps qu'un encadrement de $\mu_m$:
$2^{m/2}\le\mu _m \le 2^{(m+1)/2}$.
En 1983, Patterson et Wiedemann \cite{pw} ont montr\'e que l'on peut faire
mieux pour 
${\cal R}_{{15}}$ en exhibant une fonction bool\'eenne telle que
$\mu_m\le{27\over32}\sqrt 2\, 2^{15/2}$. Ils ont conjectur\'e que 
$\mu_m\sim 2^{m/2}$.
Rempla\c cons la fonction bool\'eenne $g$ par son exponentielle
$$f(x)=\left\{\begin{array}{rl}
1&\biq{si}g(x)=0\cr
-1&\biq{si}g(x)=1.\cr
\end{array}\right.$$
On  d\'efinit la transform\'ee de Fourier de $f$ par
$\ch f(\chi)=\sum_{V_m}f(x)\chi(x)$ o\`u $\chi$ est un caract\`ere de $V_m$,
\ie ici un homomorphisme de $V_m$ dans $\pm1$, de telle sorte que
  $  S(g)  = \|\ch f\|_\infty$.

 La conjecture de Patterson et Wiedemann se r\'e\'ecrit alors

\begin{conj}
\label{conj1}
Si $f$ d\'ecrit l'espace des fonctions de $V_m$ dans $\{\pm1\}$, on a
$$\lim_m\inf_f
{\|\ch f\|_\infty\over 2^{m/2}}=1.$$
\end{conj}

\subsubsection{Cas des tores sur $\R$}

Ce probl\`eme a un analogue avec les s\'eries de Fourier sur le tore (\ie sur
le groupe des nombres complexes de module \'egal \`a 1). Rempla\c cons en effet
les fonctions
$x\lfc (-1)^{v\cdot x}$ pour
$v\in V_m$, qui sont des caract\`eres de $V_m$ par des caract\`eres du tore
de la forme
$x\lfc e^{isx}$ pour
$s\in\Z$.

La conjecture peut se r\'e\'ecrire
$$\lim_n\inf {\|\sum_{0}^n a_{s,n} e^{isx}\|_\infty\over \sqrt n}=1$$
o\`u $a_{s,n}=\pm1$.
Autrement dit,
il existerait une suite de
polyn\^omes $P_n(z)$ et une suite de nombres positifs  $\eps_ n$ tendant
vers z\'ero tel que pour tout $| z | = 1$, $| P_n(z) |\le (1 +\eps_n)\sqrt n$,
o\`u $P_n(z) =\sum _{s=0}^n a_{s,n}z^s$ et $ a_{s,n}  = \pm1$.

Ce  probl\`eme a \'et\'e pos\'e par divers
 auteurs comme J. E. Littlewood \cite{lw}, et P. Erd\"os \cite{er}
qui a conjectur\'e qu'au contraire
il existe $\delta>1$ tel que quel que soient l'entier $n$ et le nombre
complexe
$z$ de module 1, on ait
$| P_n(z) |\ge \delta\sqrt n$.
 Kahane (\cite{ka}) a r\'esolu le probl\`eme pour des coefficients complexes
$ a_{s,n}$ de module 1, mais rien n'a \'et\'e fait pour le probl\`eme initial.
De plus, Kahane utilise pour r\'esoudre ce probl\`eme des exponentielles
de la forme
$e^{\pi i n^2/a}$, donc des exponentielles de formes quadratiques en $n$,
mais  dans notre cas elles ne donnent pas de r\'esultat complet pour les
dimensions $m$ impaires. Il fabrique avec cela un polyn\^ome qui r\'esout
presque le probl\`eme. Il ajuste ensuite ce polyn\^ome en utilisant un
argument de probabilit\'e.

\section{L'espace des fonctions bool\'eennes \`a une infinit\'e de variables}

Pour \'etudier asymptotiquement les fonctions bool\'eennes, on aura besoin de
la notion de fonction bool\'eenne \`a une infinit\'e de variable.

On rappelle que
$V_m=\f2^m$.
On d\'efinit une application de transition entre $V_m$ et $V_{m+1}$ par
$$\begin{array}{rcl}
\phi_m:&V_m&\lfd V_{m+1}\cr
&(x_1,\dots,x_{m})&\lfc
(x_1,\dots,x_m,0).
\end{array}$$
On d\'efinit
$V_\infty$ comme \'etant la limite inductive des $V_m$
suivant ces applications.

Donc
$V_\infty$ est isomorphe \`a $\f2^{(\N)}$, l'espace des suites infinies
d'\'el\'ements de
$\f2$ presque tous nuls. 
\bb

\subsection{L'espace $\Om$}

On d\'efinit $\Om_m$ comme \'etant l'ensemble des fonctions de $V_m$ dans
$\{\pm1\}$. 
Un \'el\'ement de $\Om_m$ est (l'exponentielle d') une fonction
bool\'eenne sur
$\f2^m$:
si $f$ et $g$ sont dans $\Om_m$, $fg\in\Om_m$.

On d\'efinit de mani\`ere duale aux $\phi_m$ des applications de transition
\begin{eqnarray*}
\Om_{m+1}&\lfd &\Om_m\\
 f&\lfc &f|_{V_m}
\end{eqnarray*}
o\`u $f|_{V_m}$ est la restriction de $f$ \`a $V_m$, 
 $$f|_{V_m}:(x_1,\dots,x_m)\lfc f((x_1,\dots,x_m,0).$$

Cette application permet de d\'efinir la limite projective
$$\Om=\Om_\infty=\lim proj\ \Om_n\simeq\{\pm1\}^{\f2^{(\N)}}$$ 
et les applications
$\pi_n:\Om_{\infty}\lfd \Om_n: f\lfc f|_{V_m}$.
 
On munit cet espace d'une topologie telle que
les $\pi_n^{-1}({\bf 1})$ forment un syst\`eme fondamental de voisinages de
l'origine o\`u ${\bf 1}$ est la fonction donnant \`a tous les points de $V_m$
l'image 1. Il est alors compact.

\bb

\subsection{L'espace des probabilit\'e $\Om$}

On peut munir l'espace $\Om$ d'une structure de probabilit\'e.

On d\'efinit une tribu ${\cal A}_m$ sur $\Om_m$ en prenant pour
${\cal A}_m$ l'ensemble des parties ${\cal P}(\Om_m)$ de $\Om_m$.
L'espace $\Om_m$ est muni de la probabilit\'e uniforme.
\bb
On d\'efinit la tribu ${\cal A}$ sur $\Om$ en prenant pour
${\cal A}$ la $\sigma$-alg\`ebre engendr\'ee par $\bigcup {\cal A}_m$.
On peut d\'efinir une probabilit\'e sur  cet espace $\Om$.  Pour chaque $f\in
\Om_m$, la probabilit\'e de l'\'ev\'enement $\pi_m^{-1}f$ est donn\'ee par
$\sg P (\pi_m^{-1}f)={1\over 2^q}$ o\`u $q=|V_m|=2^m$.

On notera $\E(X)$ l'esp\'erance d'une variable al\'eatoire $X$ sur $\Om$ ou
sur $\Om_m$:
$$\E(X)=\int_\Om Xd\sg P.$$

\bb

\subsection{Transformation de Fourier}

Notons $\ch V_m$ (resp. $\ch V_\infty$)
l'ensemble des caract\`ere de $V_m$ (resp. $V_\infty$).
Le groupe $V_\infty$, muni de la topologie discr\`ete, est en dualit\'e avec le
groupe
$\ch V_\infty$ qui est compact et totalement discontinu.

La transformation de Fourier est d\'efinie sur les fonctions sur $V_m$ \`a
valeurs complexes: \`a une fonction $f$ de $V_m$ dans $\C$, elle fait
correspondre une fonction
$\ch f$ de
$\ch V_m$ dans $\C$ par
$$\ch f(\chi)=\sum_{x\in V_m}f(x)\chi(x)$$
si $\chi$ est dans $\ch V_m$.

Elle se prolonge aux fonctions sur le groupe 
$V_\infty$ \`a valeurs complexes et  transforme ces fonctions en
distributions  \`a valeurs complexes sur le groupe dual
$\ch V_\infty$. Une distribution sur l'espace  $\ch V_\infty$ est une
forme lin\'eaire sur l'espace des fonctions complexes localement
constantes sur
$\ch V_\infty$ (cf. Bruhat, \cite{br}). Si
$\rho$ est une fonction test
\ie une fonction sur
$V_\infty$ qui ne prend qu'un nombre fini de valeurs non nulles,
on a
$$\sum_{x\in V_\infty}f(x)\rho(x)=\int_{\ch V_\infty}\ch f(\chi)
\ch\rho(\chi)d\chi$$ 
o\`u $d\chi$ est la mesure de Haar de $\ch V_\infty$, de masse 1.
L'expression pr\'ec\'edente a un sens si l'on remplace, comme on peut le
faire,
$V_\infty$ (et
$\ch V_\infty$) par
$V_m$ (et $\ch V_m$) pour $m$ assez grand.

\section{Etude de $\|\ch f\|_\infty$}

La relation de Parseval donne
$$q=\sum_{x\in V_m} f(x)^2=\int \ch f(\chi)^2 d\chi\le \|\ch
f\|_\infty^2$$
donc
$\|\ch f\|_\infty$ est sup\'erieur \`a $\sqrt q$.
Il est au plus \'egal \`a $q$ car
$$|\ch f(\chi)|=\Big|\sum_{x\in V_m} f(x)\chi(x)\Big|\le q.$$ 
On va montrer qu'en fait $\|\ch f\|_\infty$ est souvent voisin de $\sqrt
q$.
On montre d'abord le lemme suivant.
\begin{lem}
\label{lprep}
Si $f$ d\'esigne une fonction de $V_m$ \`a valeurs dans $\{\pm1\}$, $\chi$ un
caract\`ere de $V_m$, et $\lambda$ un r\'eel on a
$$e^{{\lambda^2q\over 2}-{\lambda^4q}}\le
{\cal E}(e^{\lambda\ch
f(\chi)})\le e^{q\lambda^2/2}.
 $$

\end{lem}

{\sl D\'emonstration} --

En effet, l'exponentielle s'\'ecrit comme un produit:
\begin{eqnarray*}
{\cal E}(e^{\lambda\ch f(\chi)})
&=&{\cal E}(e^{\lambda\sum_{x\in V_m}f(x)\chi(x)})
={\cal E}(\prod_{x\in V_m}e^{\lambda f(x)\chi(x)}).\cr
\end{eqnarray*}
Ecrivons que les variables al\'eatoires $e^{\lambda f(x)\chi(x)} $ sont
ind\'ependantes:
\begin{eqnarray*}
{\cal E}(\prod_{x\in V_m}e^{\lambda
f(x)\chi(x)})
&=&\prod_{x\in{V_m}}{\cal E}(e^{\lambda f(x)\chi(x)}).
\end{eqnarray*}
On v\'erifie que, pour $x$ fix\'e, on a
$$
{\cal E}(e^{\lambda f(x)\chi(x)})
=\cosh({\lambda })
.$$
Comme
$1+u>e^{u-{1\over2}u^2}$
si $u>0$
on a
\begin{equation}
\label{cosh}
e^{{\lambda^2\over 2}-{\lambda^4\over 8}}\le1+{\lambda^2\over 2}\le\cosh
\lambda\le e^{\lambda^2/2}.
\end{equation}
d'o\`u
$$e^{{\lambda^2q\over 2}-{\lambda^4q}}\le
e^{{\lambda^2q\over
2}-{\lambda^4q\over 8}}=
(e^{{\lambda^2\over 2}-{\lambda^4\over 8}})^q\le{\cal E}(e^{\lambda\ch
f(\chi)})\le e^{q\lambda^2/2}. $$

\subsection{Majoration de $\|\ch f\|_\infty$}

Une variante du th\'eor\`eme 1 p. 68 du livre de Kahane \cite{ka} donne le
r\'esultat suivant.

\begin{thm}
\label{th1}
Si $f$ est une fonction de $V_m$ dans $\{\pm1\}$, et $\kappa$ un r\'eel
positif, on a
$$\sg P\pt{{{\|\ch
f\|_\infty\ge \pt{2q(\kappa + \log
 q)}^{1/2}}}}
\le 2{e^{-\kappa}}.$$
\end{thm}

{\sl D\'emonstration} --

Remarquons que
$\|\ch f\|_\infty= \ch f(\chi)\hbox{ ou }-\ch f(\chi)$ pour au moins un
valeur de 
$\chi$. Donc
$$e^{\lambda\|\ch f\|_\infty}\le e^{\lambda\ch f(\chi)}+e^{-\lambda\ch
f(\chi)}$$
 pour au moins un
valeur de 
$\chi$ et par cons\'equent
 $$e^{\lambda\|\ch f\|_\infty}\le q\int_{\ch{V_m}} (e^{\lambda\ch
f(\chi)}+e^{-\lambda\ch f(\chi)})d\chi$$
d'o\`u, 
$$\displaylines{
{\cal E}\pt{e^{\lambda\|\ch f\|_\infty}}
\le q{\cal E}\pt{\int_{\ch{V_m}}
(e^{\lambda\ch f(\chi)}+e^{-\lambda\ch f(\chi)})d\chi}
\hfill\cr\hfill
\le q\pt{\int_{\ch{V_m}}
({\cal E}(e^{\lambda\ch f(\chi)})+{\cal E}(e^{-\lambda\ch
f(\chi)}))d\chi}\cr 
}$$
par inversion des sommations.
D'apr\`es le lemme \ref{lprep}
$${\cal E}\pt{e^{\lambda\|\ch f\|_\infty}}
\le 2q\pt{\int_{\ch{V_m}}
e^{q\lambda^2/2}d\chi}
\le 2qe^{q\lambda^2/2}.$$
On a donc
$${\cal E}\pt{{e^{\lambda\|\ch f\|_\infty}\over 2qe^{q\lambda^2/2}}}
\le 1$$
ou
$${\cal E}\pt{\exp\pt{{\lambda\|\ch f\|_\infty-{q\lambda^2/2}-\log(2q)}}}
\le 1.$$
Multiplions chaque membre par $2e^{-\kappa}$ o\`u $\kappa$ est un r\'eel
positif:
$${\cal E}\pt{\exp\pt{{\lambda\|\ch
f\|_\infty-{q\lambda^2/2}-\log q{-\kappa}}}}
\le 2{e^{-\kappa}}.$$
Par cons\'equent
$$\sg P\bigg(\exp\pt{{\lambda\|\ch
f\|_\infty-{q\lambda^2/2}-\log q{-\kappa}}}\ge 1\bigg)
\le 2{e^{-\kappa}}$$
\ie
$$\sg P\pt{{{\|\ch
f\|_\infty\ge{q\lambda/2}+{\log q+\kappa\over \lambda}}}}
\le 2{e^{-\kappa}}.$$
Choisissons
$\lambda=\pt{{2\kappa + 2\log q\over q}}^{1/2}$.
Cela donne le r\'esultat du th\'eor\`eme.

\begin{cor}
 On a presque-s\^urement
$$\limsup_q {\|\ch {\pi_m (f)}\|_\infty\over
2^{m/2}\sqrt{m}}\le\sqrt{2\log 2}$$ o\`u $f$ d\'ecrit  l'espace
 $\Om$.
\end{cor}

{\sl D\'emonstration} --

On prend $\kappa=\eta\log (q)$ avec $\eta>0$ dans le th\'eor\`eme pr\'ec\'edent.
On obtient
$$\sg P\pt{{{\|\ch
{\pi_m f}\|_\infty\ge \Big(2q(\eta+1)\log q \Big)^{1/2}}}}
\le {2\over q^\eta}.$$
La somme pour $m\in\N$ s'\'ecrit donc:
$$\sum_{m\in\N} \sg P\pt{{{\|\ch
{\pi_m f}\|_\infty\ge \Big(2q(\eta+1)\log q \Big)^{1/2}}}}<\infty.$$
D'apr\`es le lemme de Borel-Cantelli (cf. Kahane, \cite{ka}, \S\ 1.6, par
exemple), on en d\'eduit que, presque-s\^urement, pour
$q$ assez grand on a
$${{\|\ch
{\pi_m f}\|_\infty< \Big(2q(\eta+1)\log q \Big)^{1/2}}}.$$
Cette assertion \'etant valable pour tout $\eta$ plus grand que 0, on
peut faire tendre $\eta$ vers 0 et on a
presque-s\^urement
pour
$q$ assez grand 
$${\|\ch
{\pi_m f}\|_\infty\le (2q\log
q)}^{1/2}$$
donc presque-s\^urement
$$\limsup_m{\|\ch {\pi_m f}\|_\infty\over \sqrt{q\log (q)}}\le \sqrt 2.$$

\begin{rem}
En particulier, pour $m$ donn\'e les fonctions bool\'eennes sont en majorit\'e 
d'amplitude spectrale inf\'erieure \`a $\sqrt{2q\log
(q)}=2^{m+1\over2}\sqrt{m\log (2)}$ \`a $o(1)$ pr\`es. Carlet, et d'autre part
Olej\'ar et Stanek obtiennent le r\'esultat du th\'eor\`eme \ref{th1} \`a l'aide
d'ap\-proxi\-mations de sommes de coefficients binomiaux
\cite{ca2, os}. 
\end{rem}

\subsection{Minoration de $\|\ch f\|_\infty$}

On aura besoin des lemmes suivants.

\begin{lem}
\label{lcar}
Si $f$ d\'esigne une fonction de $V_m$ \`a valeurs dans $\{\pm1\}$,
 $\chi$ et $\xi$ deux
caract\`eres de $V_m$, et $\lambda$ un r\'eel, les majorations suivantes sont
r\'ealis\'ees:
$${\cal E}(e^{\lambda(\ch f(\chi)+\ch f(\xi))})\le 
\left\{
\begin{array}{rl}
e^{\lambda^2q}&\biq{si}\chi\ne\xi\cr
e^{2\lambda^2q}&\biq{si}\chi=\xi\ .\cr
\end{array}
\right.$$
\end{lem}

{\sl D\'emonstration} --

La d\'efinition de la transformation de Fourier, permet d'\'ecrire:
\begin{eqnarray*}
{\cal E}(e^{\lambda(\ch f(\chi)+\ch f(\xi))})
&=&{\cal E}(e^{\lambda\sum_{x\in\f{2}^m}f(x)(\chi(x)+\xi(x))})\cr
&=&{\cal E}(\prod_{x\in\f{2}^m}e^{\lambda f(x)(\chi(x)+\xi(x))})
=\prod_{x\in\f{2}^m}{\cal E}(e^{\lambda f(x)(\chi(x)+\xi(x))})
\end{eqnarray*}
puisque les variables al\'eatoire $f(x)(\chi(x)+\xi(x))$ sont ind\'ependantes
pour $x\in V_m$.
D'apr\`es le lemme \ref{lprep}, le dernier produit est \'egal \`a 
\begin{eqnarray*}
\prod_{x\in\f{2}^m}\cosh({\lambda f(x)(\chi(x)+\xi(x))})
&=&\prod_{x\in\f{2}^m}\cosh({\lambda (\chi(x)+\xi(x))})\cr
&\le&\prod_{x\in\f{2}^m}e^{\lambda^2 (\chi(x)+\xi(x))^2/2}
\end{eqnarray*}
d'apr\`es la relation (\ref{cosh}).
Ce dernier terme est \'egal \`a 
\begin{eqnarray*}
\prod_{x\in\f{2}^m}e^{\lambda^2 (1+\chi\xi(x))}
&=&\prod_{x\in\f{2}^m}e^{\lambda^2}e^{\lambda^2 \chi\xi(x)}
=e^{\lambda^2q}\prod_{x\in\f{2}^m}e^{\lambda^2 \chi\xi(x)}\cr
&=&e^{\lambda^2q}e^{ \lambda^2\sum_{x\in\f{2}^m}\chi\xi(x)}
=\left\{
\begin{array}{rl}
e^{\lambda^2q}&\biq{si}\chi\ne\xi\cr
e^{2\lambda^2q}&\biq{si}\chi=\xi\cr
\end{array}
\right.
\end{eqnarray*}
en utilisant l'annulation de la somme des valeurs des caract\`eres non
triviaux.

On aura \'egalement besoin d'une in\'egalit\'e \'el\'ementaire:
\begin{lem}
\label{ineg}
Si $X$ est une variable al\'eatoire de carr\'e int\'egrable et si
$0<\lambda<1$, on a
$$\sg P\Big(X\ge\lambda\E(X)\Big)\ge(1-\lambda)^2{\E^2(X)\over \E(X^2)}.$$
\end{lem}

{\sl D\'emonstration} --

 Voir par exemple Kahane \cite{ka}, \S\ 1.6.

\bb

Le th\'eor\`eme suivant donne une minoration de la probabilit\'e que $\|\ch
f\|_\infty$ soit assez grand. Il est inspir\'e de Salem et Zygmund
\cite{sz} qui traitent le cas du tore. Voir aussi l'article de B. Kashin
et L. Tsafriri
\cite{kt}.

\begin{thm}
\label{minf}
Si $f$ d\'esigne une fonction de $V_m$ \`a valeurs dans $\{\pm1\}$,
et si $0<\alpha<1$ et $0<\eta<1-\alpha^2$, alors il existe une constante
$B$ positive et ne d\'ependant que de $\alpha$ et $\eta$ telle que 
$$\sg P \Bigg(\|\ch f\|_\infty> \Big({\alpha \over 2}-{\eta \over
\alpha}-\alpha^3 {\log q\over q}\Big)\sqrt{q\log q}\Bigg)>1-{B\over
q^\eta}.
$$
\end{thm}

{\sl D\'emonstration} --

D\'efinissons la variable al\'eatoire
$I_q=  \int_{\ch{V_m}} \exp(\lambda \ch f(\chi))$.
Le lemme \ref{lprep} permet  de minorer $\E(I_q)$:
$$\E(I_q)=\int_{\ch{V_m}} \E(\exp(\lambda \ch f(\chi)))
\ge\int_{\ch{V_m}} e^{{\lambda^2q\over 2}-{\lambda^4q}}
= e^{{\lambda^2q\over 2}-{\lambda^4q}}.
$$
De plus
$$I_q(\chi)^2=  \int_{\ch{V_m}} \exp(\lambda \ch f(\chi))
\int_\xi \exp(\lambda \ch f(\xi))
=  \int_{\chi,\xi }\exp(\lambda (\ch f(\chi)+\ch f(\xi)))
$$
d'o\`u, d'apr\`es le lemme \ref{lcar} pr\'ec\'edent
$$\displaylines{
\E(I_q(\chi)^2)
=  \int_{\chi,\xi }\E(\exp(\lambda (\ch f(\chi)+\ch f(\xi))))
\hfill\cr\hfill
\le  \int_{\chi,\xi }\exp(q\lambda^2)
+{1\over q}\int_{\chi }\exp(2q\lambda^2)
=\Big(1+{\exp(q\lambda^2)\over q}\Big)\exp(q\lambda^2).
}
$$
Donc, d'apr\`es l'in\'egalit\'e du lemme \ref{ineg}, si $\eta$ est un r\'eel
positif
$$\displaylines{
\sg P \Big(I_q>q^{-\eta} e^{{\lambda^2q\over
2}-{\lambda^4q}}\Big) 
\ge(1-q^{-\eta})^2{e^{{\lambda^2q}-{2\lambda^4q}}\over 
\Big(1+{\exp(q\lambda^2)\over q}\Big)\exp(q\lambda^2)}
\hfill\cr\hfill
\ge(1-2q^{-\eta}){e^{-{2\lambda^4q}}}\Big(1-{\exp(q\lambda^2)\over
q}\Big)}$$ 
si
${\exp(q\lambda^2)\over q}<1$.
Cette condition est satisfaite si on choisit 
$$\lambda=\alpha\pt{\log q\over q}^{1/2}$$
avec
$0<\alpha<1$.
Ce choix de $\lambda$ permet de calculer $e^{-{2\lambda^4q}}$:
$$2\lambda^4q=2\alpha^4\pt{\log q\over q}^{2}q=2\alpha^4{\pt{\log
q}^{2}\over q} <2{\pt{\log q}^{2}\over q}$$
d'o\`u
$$\sg P \Big(I_q>q^{-\eta} e^{{\lambda^2q\over
2}-{\lambda^4q}}\Big) 
\ge\Big(1-{2\over q^\eta}\Big)\Big(1-2{\pt{\log
q}^{2}\over q}\Big)\Big(1-q^{\alpha^2-1}\Big)>\Big(1-{B\over
q^\eta}\Big)$$ 
pour une certaine constante $B$
si $\eta<1-\alpha^2$.

Il est \'evident que
$$\exp(\lambda \|\ch f\|_\infty)\ge 
\int_{\ch{V_m}}\exp(\lambda \ch f(\chi)).$$
D'o\`u
$$\sg P \Big(\exp(\lambda \|\ch f\|_\infty)>q^{-\eta} e^{{\lambda^2q\over
2}-{\lambda^4q}}\Big) 
\ge
\sg P \Big(I_q>q^{-\eta} e^{{\lambda^2q\over
2}-{\lambda^4q}}\Big) .$$ 
L'\'ev\'enement du membre de gauche peut encore s'\'ecrire
$$ \|\ch f\|_\infty> {{\lambda q\over 2}-{\lambda^3q}-\eta\log
q/\lambda}.$$ 
Ou encore
$$ \|\ch f\|_\infty> {{\alpha \over 2}\sqrt{q\log
q}-{\lambda^3q}-\eta\log
q/\lambda}.$$
Majorons le deuxi\`eme terme de cette somme:
$$\lambda^3 q=\alpha^3 \pt{\log q\over q}^{3/2}q
=\alpha^3 {(\log q)^{3/2}\over q^{1/2}}
=\alpha^3 {\log q\over q}\sqrt{q\log q}.
$$
Enfin le troisi\`eme terme vaut 
$$\eta\log q/\lambda={\eta\log q\over \alpha}\pt{q\over \log 
q}^{1/2} ={\eta \over \alpha}\sqrt{q\log q}.$$
L'\'ev\'enement en question s'\'ecrit donc
 \begin{eqnarray*}
\|\ch f\|_\infty> \Big({\alpha \over 2}-{\eta \over \alpha}-\alpha^3 {\log
q\over q}\Big)\sqrt{q\log q}
\end{eqnarray*}
ce qui termine la d\'emonstration.

\begin{cor}
 On a presque-s\^urement
$$\liminf_m{\|\ch {\pi_m(f)}\|_\infty\over2^{m/2}\sqrt{m}}\ge {\log 2
\over 2}$$
o\`u $f$ d\'ecrit les \'el\'ements de l'espace
 $\Om$.
\end{cor}

En effet, en faisant la somme pour $m\in\N$ des in\'egalit\'es donn\'ees par le
th\'eor\`eme pr\'ec\'edent, on obtient
$$\sum_m\sg P \bigg(\|\ch {\pi_m(f)}\|_\infty< \Big({\alpha \over 2}-{\eta
\over
\alpha}-\alpha^3 {\log q\over q}\Big)\sqrt{q\log q}\bigg)<
\sum_m {B\over q^\eta}<\infty.
$$
Donc,  le lemme de Borel-Cantelli nous dit que
p.s. 
$$\|\ch {\pi_m(f)}\|_\infty> \Big({\alpha \over 2}-{\eta \over
\alpha}-\alpha^3 {\log q\over q}\Big)\sqrt{q\log q}$$ sauf pour un nombre
fini de $q$, 
\ie p.s.
$$\liminf_m{\|\ch {\pi_m(f)}\|_\infty\over\sqrt{q\log q}}> {\alpha
\over 2}-{\eta
\over
\alpha}.$$
On peut faire tendre $\alpha$ vers 1 et $\eta $ vers 0.
On obtient
$$p.s.\qquad\liminf_m{\|\ch {\pi_m(f)}\|_\infty\over\sqrt{q\log q}}\ge {1
\over 2}.$$

\section{Etude de $\|\ch f\|_4$}

Reprenons  l'id\'ee de
D. Newman et J. Byrnes \cite{nb}. Ils ont remarqu\'e que, dans le cas  des
s\'eries de Fourier  sur $\Z$, la norme dans $L^4$ de $\sum_n\pm e^{int}$
avait une expression agr\'eable.
Il en va de m\^eme de $\|\ch f\|_4$ pour $f:V_m\fd\{\pm1\}$.
On remarque que 
\begin{equation}
\label{ineg2}
\|\ch f\|_2\le\|\ch f\|_4\le\|\ch f\|_\infty.
\end{equation}
En effet,
la premi\`ere in\'egalit\'e   vient de ce que   les fonction $\ch f$ sont
d\'efinies sur un espace de mesure \'egale \`a 1.  
Par cons\'equent, la conjecture \ref{conj1}
 implique une conjecture plus faible:
\begin{conj}
Si $f$ d\'ecrit l'espace des fonctions de $V_m$ dans $\{\pm1\}$, on a
$$\lim_m\inf_{f\in V_m} {\|\ch f\|_4\over 2^{m/2}}=1.$$
\end{conj}

L'id\'ee d'\'etudier $\|\ch f\|_4$ n'est pas nouvelle puisque C. Carlet a
propos\'e d'\'etudier la non-lin\'earit\'e des fonctions bool\'eennes  par
les  moments d'ordre sup\'erieur de leur transform\'ees de Fourier \cite{ca}.
Cela a \'egalement \'et\'e \'etudi\'e par Xian-Mo Zhang et Yuliang Zheng  \cite{zz},
ou par P. St\u anic\u a \cite{st}  sous le nom de 
``somme des carr\'es''.  
On peut voir \'egalement l'article de P.~Langevin et P.~Sol\'e \cite{ls} qui
appliquent cette notion \`a une cubique.

\subsection{L'expression de $\|\ch f\|_4$}

On obtient l'expression simple suivante pour $\|\ch f\|_4$.

\begin{lem}
\label{som}
Si  $f$ est une fonction de $V_m$ \`a valeurs dans $\pm1$,
$$\|\ch f\|_4^4 = \sum_{x_1+x_2+x_3+x_4=0}
f(x_1)f(x_2)f(x_3)f(x_4).$$
\end{lem}

{\sl D\'emonstration.} --

D\'ecomposons $\ch f^4$ et  inversons l'ordre de la somme et de
l'int\'egrale:
\begin{eqnarray*}
\|\ch f\|_4^4 &=& \int_{\ch V_m} \ch f^4 d\chi\cr
&=& \int_{\ch V_m} \pt{\sum_{{ V_m}}f(x_1)\chi(x_1)}
\pt{\sum_{{ V_m}}f(x_2)\chi(x_2)}
\cr&&\hskip4cm
\pt{\sum_{{ V_m}}f(x_3)\chi(x_3)}
\pt{\sum_{{ V_m}}f(x_4)\chi(x_4)}d\chi\cr
&=& \sum_{x_1,x_2,x_3,x_4}
f(x_1)f(x_2)f(x_3)f(x_4)\int_{\ch V_m}\chi(x_1+x_2+x_3+x_4)d\chi\cr
&=& \sum_{x_1+x_2+x_3+x_4=0}
f(x_1)f(x_2)f(x_3)f(x_4).\cr
\end{eqnarray*}

\subsubsection{R\'e\'ecriture de la conjecture \ref{conj1}}

D\'ecomposons la somme donn\'ee dans le lemme \ref{som}:
$$\displaylines{
\sum_{x_1+x_2+x_3+x_4=0}
f(x_1)f(x_2)f(x_3)f(x_4)
\hfill\cr\hfill
=q^2+\sum_{a\ne0}\ \sum_{x_1+x_2=x_3+x_4=a}
f(x_1)f(x_2)f(x_3)f(x_4)}$$
par suite
$$\displaylines{
\sum_{x_1+x_2=x_3+x_4=a}f(x_1)f(x_2)f(x_3)f(x_4)\\
\hfill\cr\hfill
=\sum_{x_1+x_2=a}f(x_1)f(x_2)\sum_{x_3+x_4=a}f(x_3)f(x_4)
=\pt{\sum_{x_1+x_2=a}f(x_1)f(x_2)}^2\ge0.
}$$
D\'efinissons les variables al\'eatoires \`a valeurs dans $\C$ et d\'ependant de
$a$ dans
$V_m$:
\begin{equation}
\label{xa}
\dps X_a=\pt{\sum_{x_1+x_2=a}f(x_1)f(x_2)}^2.
\end{equation}
D'o\`u
$$
\|\ch f\|_4^4
-q^2={\sum_{\textstyle{a\ne0\atop a\in V_m}}X_a}.$$
La conjecture \ref{conj1} se r\'e\'ecrit de la mani\`ere suivante.

\begin{conj}
Pour tout $\eps>0$, il existe $q$ non carr\'e et $f$ dans $V_m$ tels que
$\dps\sum_{a\ne0}X_a<\eps q^2$ o\`u $X_a$ est donn\'e par la formule
(\ref{xa}).
\end{conj}

On comparera avec le probl\`eme 14.3 dans l'article de Tamas Erd\'elyi
\cite{er0}.

\subsection{Calculs d'esp\'erances}

Remarquons que $\|\ch f\|_4^4$ est compris entre $q^2$ et $q^3$. En effet,
la premi\`ere in\'egalit\'e  $\|\ch f\|_4^4\ge q^2$ vient de l'in\'equation
\ref{ineg2}.   Le lemme \ref{som} implique d'autre part que
 $\|\ch f\|_4^4\le q^3$. Cf. \cite{zz}, th\'eor\`eme 6 ou \cite{cccf}, th\'eor\`eme
1.

On peut \'egalement d\'eduire du lemme \ref{som} le calcul de
$\E(\|\ch f\|_4^4)$
et de
$\E(\|\ch f\|_4^8)$.
On utilise pour cela les lemmes suivants. Remarquons d'abord que
$${\cal E}(\|\ch f\|_4^4 )= \sum_{x_1+x_2+x_3+x_4=0}{\cal
E}(f(x_1)f(x_2)f(x_3)f(x_4)) .$$

\begin{lem}
\label{elim}
On a
$$
\E(f(x_1)f(x_2)\dots
f(x_r))=\E(f(x_3)f(x_4)\dots
f(x_r))$$
si $x_1=x_2$.
\end{lem}

{\sl D\'emonstration} --

Si $x_1=x_2$, on a $
f(x_1)f(x_2)\dots
f(x_r)=f(x_3)f(x_4)\dots
f(x_r)$.

\begin{lem}
\label{pair}
On a 
$$
\E(f(x_1)f(x_2)f(x_3)\dots
f(x_r))=0\biq{ou}1. $$
L'esp\'erance
$\E(f(x_1)f(x_2)f(x_3)\dots
f(x_r))$ est \'egale \`a 1 \ssi pour chaque $y\in \f2^n$ l'ensemble des $x_i$
\'egaux \`a
$y$ a un cardinal pair, \ie
\ssi il existe une partition de $\{x_1,x_2,x_3,\dots,x_r\}$ en couples
form\'es d'\'el\'ements \'egaux.
\end{lem}

{\sl D\'emonstration} ---

L'application successive du lemme pr\'ec\'edent permet de se ramener au cas
o\`u tous les $x_i$ sont distincts.
Les variables al\'eatoire $f(x_i)$ sont alors ind\'ependantes, donc
$$\E(f(x_1)f(x_2)f(x_3)\dots
f(x_r))=\E(f(x_1))\E(f(x_2))\E(f(x_3))\dots
\E(f(x_r)).$$
De plus
$f(x_1)=1\hbox{ ou } -1$ avec la probabilit\'e $1/2$. D'o\`u $\E(f(x_1))=0$.
Donc si tous les $x_i$ sont distincts,
$\E(f(x_1)f(x_2)f(x_3)\dots f(x_r))=0$ sauf si la suite des $x_i$ est
vide auquel cas elle vaut $\E(1)=1$.

\bb

Posons
$$E(a_1,\dots,a_r)=\sum_{x_1,\dots,x_r}\E\Big(f(x_1)f(x_1+a_1)\dots
f(x_r)f(x_r+a_r)\Big).$$

\begin{lem}
\label{red}
Si $a$ et les $a_i$ sont dans $V_m$, on a, pour $a\ne0$
$$E(a,a_1,\dots,a_r)\le 2\sum_{1\le i\le r}
E(a_1,\dots,a_i+a,\dots,a_r).$$
\end{lem}

{\sl D\'emonstration} --

Appliquons  le lemme \ref{pair}:
\begin{eqnarray*}
&&E(a,a_1,\dots,a_r)\\&=
&\sum_{x,x_1,\dots,x_r}\E\Big(f(x)f(x+a)f(x_1)f(x_1+a_1)\dots
f(x_r)f(x_r+a_r)\Big)
\cr 
&=&\sum_{x_1,x_2,\dots,x_r}
\quad\sum
\E(f(x)f(x+a)f(x_1)f(x_1+a_1)
\dots f(x_r)f(x_r+a_r))
\end{eqnarray*}
o\`u la derni\`ere somme est sur les $x$ dans
$\{x_1,x_1+a_1,\dots,x_r,x_r+a_r\}$.
Si $x=x_1$,  le lemme \ref{elim} permet d'\'ecrire
\begin{eqnarray*}
&&\hskip-2cm\E(f(x)f(x+a)f(x_1)f(x_1+a_1)\dots
f(x_r)f(x_r+a_r))\\
&=&\E(f(x_1)f(x_1+a)f(x_1)f(x_1+a_1)\dots
f(x_r)f(x_r+a_r))\cr
&=&\E(f(x_1+a)f(x_1+a_1)\dots
f(x_r)f(x_r+a_r))\cr
&=&\E(f(t)f(t+a_1+a)\dots
f(x_r)f(x_r+a_r))
\end{eqnarray*}
en posant $t=x_1+a_1$.
Si $x=x_1+a_1$, on a de m\^eme
$$\displaylines{
\E(f(x)f(x+a)f(x_1)f(x_1+a_1)\dots
f(x_r)f(x_r+a_r))\hfil\cr\hfil
=\E(f(x_1)f(x_1+a+a_1)f(x_2)f(x_2+a_2)\dots
f(x_r)f(x_r+a_r))
}$$
d'o\`u le lemme.

\subsubsection{Les esp\'erances de $X_a$, $X_a^2$, $X_aX_b$}

\begin{prop}
Si $a$ est un \'el\'ement non nul de $V_m$, on a
${\cal E}(X_a)=2q.$
\label{esp1}
\end{prop}
{\sl D\'emonstration} ---

D\'ecomposons ${\cal E}(X_a)$:
\begin{eqnarray*}
{\cal E}\pt{\sum_{x_1+x_2=a}f(x_1)f(x_2)}^2
&=&
{\cal E}\pt{\sum_{x_1\in V_m}f(x_1)f(x_1+a)}^2\\
&=&\sum_{x_1,x_2}{\cal E}\pt{f(x_1)f(x_1+a)f(x_2)f(x_2+a)}.\cr
\end{eqnarray*}
Les valeurs de $(x_1,x_2)$ qui rendent 
${\cal E}\pt{f(x_1)f(x_1+a)f(x_2)f(x_2+a)}$ \'egale \`a 1 sont
$x_1=x_2$, et $x_1=x_2+a$. Il y en a $q$ dans les deux cas.
D'o\`u la proposition.

\begin{prop}
Si $a$ est un \'el\'ement non nul de $V_m$, on a
${\cal E}(X_a^2)\le12q^2.$
\end{prop}
{\sl D\'emonstration} ---

D\'ecomposons $X_aX_a$:
\begin{eqnarray*}
X_aX_a&=&\pt{\sum_{x}f(x)f(x+a)}^2
\pt{\sum_{x}f(x)f(x+a)}^2\cr
&=& \pt{\sum_{x}f(x)f(x+a)}\pt{\sum_{y}f(y)f(y+a)}\\
&&\hskip4cm\pt{\sum_{z}f(z)f(z+a)}\pt{\sum_{t}f(t)f(t+a)}\cr
&=& \sum_{x,y,z,t}f(x)f(x+a)f(y)f(y+a)
f(z)f(z+a)f(t)f(t+a) .\cr
\end{eqnarray*}
D'o\`u, d'apr\`es le lemme  \ref{red} 
\newdimen\ecart
\ecart=5cm
\begin{eqnarray*}
\E(X_a^2)
&=& \sum_{x,y,z,t}\E\pt{f(x)f(x+a)f(y)f(y+a)
f(z)f(z+a)f(t)f(t+a)} \cr
&=& E(a,a,a,a)\\
&\le& 2 E(0,a,a)+2E(a,0,a)+2E(a,a,0)\\
&=& 6q E(a,a).
\end{eqnarray*}
d'apr\`es le lemme \ref{elim}. Ce dernier terme vaut
 $12q^2 $
d'apr\`es la proposition \ref{esp1}.

\begin{prop}
Si $a$ et $b$ sont des \'el\'ements non nuls et distincts de $V_m$, on a
${\cal E}(X_aX_b)\le 4q^2+32q.$
\end{prop}
{\sl D\'emonstration} ---

D\'ecomposons $X_aX_b$:
\begin{eqnarray*}
X_aX_b
&=& \pt{\sum_{x}f(x)f(x+a)}^2
\pt{\sum_{x}f(x)f(x+b)}^2\cr
&=& \pt{\sum_{x}f(x)f(x+a)}\pt{\sum_{y}f(y)f(y+a)}\\
&&\hskip4cm\pt{\sum_{z}f(z)f(z+b)}\pt{\sum_{t}f(t)f(t+b)}\cr
&=& \sum_{x,y,z,t}f(x)f(x+a)f(y)f(y+a)
f(z)f(z+b)f(t)f(t+b) .\cr
\end{eqnarray*}

On a donc d'apr\`es de lemme \ref{pair}
\newdimen\comm
\comm=6cm
\begin{eqnarray*}
\E(X_aX_b)
&=& \sum_{x,y,z,t}\E\pt{f(x)f(x+a)f(y)f(y+a)
f(z)f(z+b)f(t)f(t+b)} \cr
&=& E(a,a,b,b) \cr
&\le& 2E(0,b,b) +2E(a,b+a,b)+2E(a,b,b+a)\cr
&=& 2E(0,b,b) +4E(a,b+a,b).
\end{eqnarray*}
en utilisant l'\'egalit\'e $E(a,b+a,b)=E(a,b,b+a)$.
En utilisant le lemme \ref{elim}, on obtient
$$\E(X_aX_b)
\le 2qE(b,b) +4E(a,b+a,b).$$
La premi\`ere somme est calcul\'ee dans da proposition \ref{esp1}.
Calculons la deuxi\`e\-me somme, en utilisant le lemme \ref{red}.
\begin{eqnarray*}
E(a,b+a,b)\le 2 E(b,b)+2E(b+a,b+a)=8q
\end{eqnarray*}
d'apr\`es la proposition \ref{esp1}.

\subsubsection{Esp\'erances de $\|\ch f\|_4^4$ et $\|\ch f\|_4^8$}
\label{esp}

\begin{prop}
Si $f$ est une fonction de $V_m$ dans $\{\pm1\}$, alors
$\E(\|\ch f\|_4^4)=3q^2-2q$.
\end{prop}

{\sl D\'emonstration} --

En effet
$$\|\ch f\|_4^4=q^2+\sum_{a\ne0}X_a$$
Donc
$$\E(\|\ch f\|_4^4)=q^2+\sum_{a\ne0}\E(X_a)=
q^2+2q(q-1).$$

\begin{prop}
Si $f$ est une fonction de $V_m$ dans $\{\pm1\}$,
$\E(\|\ch f\|_4^8)\le
64q - 100q^2 + 28q^3 + 9q^4$.
\end{prop}

{\sl D\'emonstration} --

On a
\begin{eqnarray*}
\|\ch f\|_4^8
&=&
(q^2+\sum_{a\ne0}X_a)^2\cr
&=&
q^4+2q^2\sum_{a\ne0}X_a+\sum_{a\ne0}X_a^2+\sum_{0\ne a\ne b\ne0}X_aX_b.\cr
\end{eqnarray*}
Donc
\begin{eqnarray*}
\E(\|\ch f\|_4^8)
&=&
q^4+2q^2\sum_{a\ne0}\E(X_a)+\sum_{a\ne0}\E(X_a^2)+\sum_{0\ne a\ne
b\ne0}\E(X_aX_b)\cr
&\le&
q^4+2q^2(q-1)2q+12q^2(q-1)+(q-1)(q-2)(4q^2+32q).\cr
&=&
64q - 100q^2 + 28q^3 + 9q^4\cr
\end{eqnarray*}

\subsection{In\'egalit\'es sur $\|\ch f\|_4^4$}

\begin{prop}
Si $f$ est une fonction de $V_m$ dans $\{\pm1\}$, et $t$ un r\'eel positif,
$$\sg P\Bigg(\bigg|{\|\ch
f\|_4^4\over q^2}-3+{2\over q}\bigg|\ge t\Bigg)\le{40\over t^2q}.$$
\end{prop}

{\sl D\'emonstration} --

La variance de $\|\ch f\|_4^4$ v\'erifie, d'apr\`es le paragraphe \ref{esp}
pr\'ec\'edent
$$
var(\|\ch f\|_4^4)={\cal E}(\|\ch f\|_4^8)-{\cal E}(\|\ch f\|_4^4)^2
\le64q - 104q^2 + 40q^3.
$$

Donc en appliquant l'in\'egalit\'e de
Bienaym\'e-Tchebicheff (Voir par exemple Kahane \cite{ka}, \S\ 1.6.)
$$\sg P\bigg(\Big|\|\ch f\|_4^4-{\cal E}(\|f\|_4^4)\Big|\ge u\bigg)\le 
{var(\|f\|_4^4)\over u^2}\le{64q - 104q^2 + 40q^3\over u^2}$$
pour $u>0$.
En rempla\c cant $u$ par $q^2t$, on obtient:
$$\sg P\bigg(\Big|\|\ch f\|_4^4-3q^2+2q\Big|\ge q^2t\bigg)\le 
{64q - 104q^2 + 40q^3\over q^4t^2}\le{40\over t^2q}$$
d'o\`u le r\'esultat.

\subsection{Etude asymptotique de $\|\ch{f}\|_4$}

Pour presque tout $f$ appartenant \`a $\Om$, $\|\ch{\pi_m f}\|_4\over \sqrt
q$ a une limite donn\'ee par le corollaire suivant.

\begin{cor}
Si $f\in\Om$, on a presque s\^urement
$$\lim_m{\|\ch{\pi_m f}\|_4\over 2^{m/2}}=3^{\scr 1/4}.$$
\end{cor}

{\sl D\'emonstration} --

Faisons la somme pour $m\in\N$ des in\'egalit\'es donn\'ees par la proposition
pr\'ec\'edente:
$$\sum_m\sg P\bigg(\bigg|{\|\ch{\pi_m f}\|_4^4\over q^2}-3+{2\over
q}\bigg|\ge t\bigg)\le\sum_m{40\over t^2q}<\infty$$
par cons\'equent, le lemme de Borel-Cantelli dit que, pout $t$ donn\'e,  
on a presque s\^urement
 $$\bigg|{\|\ch{\pi_m f}\|_4^4\over q^2}-3+{2\over q}\bigg|<
t$$ sauf peut-\^etre pour un nombre fini de $q$.
Par cons\'equent, on a presque s\^urement
$$\lim_m{\|\ch{\pi_m f}\|_4^4\over q^2}=3.$$

\subsection{R\'esultats asymptotiques}

\subsubsection{Convergence de la loi de la variable al\'eatoire ${1\over
q}X_a$}

On notera $$\Phi_{X}(u)=\E\pt{\exp(iuX)}$$ la fonction
caract\'eristique d'une variable al\'eatoire $X$.
\begin{prop}
La distribution de
${1\over\sqrt q}\pt{\sum_{x\in\f2^m}f(x)f(x+a)}$
converge en loi vers la distribution 
gaussienne d'esp\'erance nulle et de variance $2$ quand $q$ tend vers
l'infini.
\end{prop}

{\sl D\'emonstration} ---

Soit $H$ l'hyperplan de l'espace vectoriel $\f2^m$ orthogonal \`a $a$.
Les variables al\'eatoires $f(x)f(x+a)$ sont ind\'ependantes pour $x\in H$ et
on a
$${1\over\sqrt q}\pt{\sum_{x\in\f2^m}f(x)f(x+a)}
={2\over\sqrt q}\pt{\sum_{x\in H}f(x)f(x+a)}.$$

Le d\'eveloppement de Taylor de $\Phi_{f(x)f(x+a)}$ \`a l'origine est donn\'e
par
\begin{eqnarray*}
\Phi_{f(x)f(x+a)}(u)&=&
1+iu\E(f(x)f(x+a))\\
&&\hskip2cm-u^2\E(f(x)f(x+a)f(x)f(x+a))/2+o(u^2)\\&=&
1-u^2/2+o(u^2)
\end{eqnarray*}
dans un voisinage de l'origine.
Donc
$$\log\Phi_{f(x)f(x+a)}(u)=
-u^2/2+o(u^2)$$
dans un voisinage $\cal V$ de l'origine.

Posons $$S_q=\sum_{x\in H}f(x)f(x+a)$$
et soit $\Phi_{S_q/\sqrt q}(u)=\E\pt{\exp(iuS_q/\sqrt q}$ la fonction
caract\'eristique de $S_q/\sqrt q$.  On a
$$\Phi_{S_q/\sqrt q}(u)=\E\pt{\exp(iuS_q/\sqrt q}
=\E\pt{\exp(iS_qu/\sqrt q}=\Phi_{S_q}(u/\sqrt q).$$
Ecrivons que les variables al\'eatoires
$f(x)f(x+a)$ sont ind\'ependantes pour $x\in H$:
\begin{eqnarray*}
\Phi_{S_q}(u/\sqrt q)
&=&\prod_{x\in H}\Phi_{f(x)f(x+a)}(u/\sqrt q)
\end{eqnarray*}
d'o\`u
\begin{eqnarray*}
\log\Phi_{S_q}(u/\sqrt q)
&=&\sum_{x\in H}\log\Phi_{f(x)f(x+a)}(u/\sqrt q)\\
&=&\sum_{x\in H}\bigl(-u^2/2 q+o(u^2/q)\bigr)\\
&=&-u^2/4 +qo(u^2/q)
\end{eqnarray*}
d\`es que $q$ est assez grand pour que $u/\sqrt q$ soit dans le voisinage
$\cal V$ de 0. Par cons\'equent, quand $q$ tend vers l'infini, la fonction
caract\'eristique
$\Phi_{S_q/\sqrt q}$ tend vers $\exp(-u^2/4) $, donc vers la fonction
caract\'eristique d'une loi normale centr\'ee, de variance $1/2$. 

La distribution de
${1\over\sqrt q}\pt{\sum_{x\in\f2^m}f(x)f(x+a)}$
\'egale \`a ${2}S_q/\sqrt q$ converge donc en loi  vers une loi normale
centr\'ee, de variance $2$.

Posons $Y_a={1\over q}X_a$.
 
\begin{prop}
La distribution de
$Y_a={1\over q}X_a$
converge en loi vers la distribution 
de densit\'e
$${1\over 2\sqrt{\pi x}}{e^{- x/4}}{\bf 1}_{(x>0)}.$$
\end{prop}

{\sl D\'emonstration} ---

En effet, la variable $Y_a$ est \'egale \`a 
$\pt{{1\over\sqrt q}\pt{\sum_{x\in(\f2^m)}f(x)f(x+a)}}^2$, et le carr\'e
d'une distribution 
gaussienne d'esp\'erance nulle et de variance $2$ est une variable
al\'eatoire de densit\'e
$${1\over 2\sqrt{\pi x}}{e^{- x/4}}{\bf 1}_{(x>0)}.$$

\subsubsection{Le th\'eor\`eme de G\"artner-Ellis}
Il s'agit, comme le le th\'eor\`eme de Cramer, d'un th\'eor\`eme concernant 
une \'evaluation des grandes d\'eviations d'une variable al\'eatoire. On pourra
se r\'ef\'erer \`a  J. Bucklew \cite{bu} ou \`a
A. Dembo et  O. Zeitouni \cite{dz}.

On a
$${1\over q^2}\|\ch f\|_4^4
-1={1\over
q}{\sum_{a\ne0}Y_a}.$$
On a vu que les variables al\'eatoires $Y_a$ avaient presque la m\^eme
distribution. Si elles  \'etaient ind\'ependantes,
on pourrait r\'esoudre la conjecture \ref{conj1} en utilisant le th\'eor\`eme de
Cramer, qui donne une \'evaluation des grandes d\'eviations pour des
variables al\'eatoire ind\'ependantes et \'equidistribu\'ees. 

Mais cela n'est pas le cas ici. Le th\'eor\`eme de G\"artner-Ellis peut
peut-\^etre s'appliquer, mais il faut v\'erifier un certain nombre
d'hypoth\`eses.
 D\'efinissons 
$$\phi_q(u)={1\over q}\log\E\pt{\exp\Big(u{\sum_{a\ne0}Y_a}\Big)}.$$
Supposons que $\phi(u)=\lim_{q\fd\infty} \phi_q(u)$ existe pour tout
$u\in\R$ (prenant \'eventuellement des valeurs infinies) et soit
diff\'erentiable sur l'ensemble
$D_\phi=\{u\mid\phi(u)<\infty\}.$

Soit aussi
$$I(x)=\sup_u(ux-\phi(u))$$
et
$$\phi'(D_\phi)=\{\phi'(u)\mid u\in D_\phi\}.$$
Le th\'eor\`eme de G\"artner-Ellis implique
$$\limsup_{q\fd\infty}{1\over q}\log\sg P\pt{{1\over
q}\sum_{a\ne0}Y_a<\eps}
\le-\inf_{x<\eps}I(x).$$
Comme pour $q$ carr\'e, il existe des fonctions courbes, la probabilit\'e que 
$\|\ch f\|_4^4
=q^2$ est sup\'erieure \`a $2^q$, donc
$$\log 2\le -\inf_{x<\eps}I(x).$$

Mais
si $]0,\eps[\subset \phi'(D_\phi)$, le th\'eor\`eme de G\"artner-Ellis implique
$$\liminf_{q\fd\infty}{1\over q}\log\sg P\pt{{1\over
q}\sum_{a\ne0}Y_a<\eps}
\ge-\inf_{x<\eps}I(x).$$
On en d\'eduirait  que pour $\eps$ donn\'e, pour tout $q$ assez grand, il
existe
$f$ tel que
$${1\over q^2}\|\ch f\|_4^4
-1<\eps.$$
donc tel que
$${1\over \sqrt q}\|\ch f\|_4-1<\eps.$$


\newpage
\parindent=0mm
{\bf Fran\c cois Rodier}

\noindent
Institut de Math\'ematiques de Luminy, 

163 Avenue de Luminy, 

Case
907, 

13288 Marseille cedex 9 -- France.

\bs
\tt rodier@iml.univ-mrs.fr

\end{document}